\documentclass[a4paper,leqno]{article}
\usepackage{latexsym,amssymb,amsmath}
\usepackage{color,amsfonts}
\usepackage[english]{babel}

\newtheorem{theorem}{Theorem}[section]
\newtheorem{lemma}[theorem]{Lemma}
\newtheorem{proposition}[theorem]{Proposition}
\newtheorem{corollary}[theorem]{Corollary}

\newtheorem{remark}[theorem]{Remark}
\newenvironment{proof}[1][Proof]{\textbf{#1.} }{\ \mathbb Rule{0.5em}{0.5em}}

\newcommand{\po}{{\mathbb P}}
\newcommand{\ep}{{\varepsilon}}
\newcommand{\ms}{{\mu^{\beta,\nu}}}
\newcommand{\msu}{{\mu^{1,\nu}}}

\newcommand{\ove}{{\overline v^\varepsilon}}

\begin{document}

\title{2D hydrodynamical systems:
       invariant measures of Gaussian type}
\author{Hakima Bessaih
\footnote{University of Wyoming, Department of Mathematics, Dept. 3036, 1000
East University Avenue, Laramie WY 82071, United States, bessaih@uwyo.edu}
\and Benedetta Ferrario
\footnote{Universit\`a di Pavia, Dipartimento di Matematica, via
  Ferrata 1, 27100 Pavia, Italy,
benedetta.ferrario@unipv.it}}

\date{}
\maketitle

\begin{abstract}
Gaussian measures $\ms$ are associated to some stochastic
2D hydrodynamical systems. They are of Gibbsian type and are
constructed
by means of some invariant quantities of the system depending
on some parameter $\beta$ (related to the 2D nature of the fluid)
and the viscosity $\nu$.
We prove the existence and the uniqueness of the global flow for the
stochastic viscous system; moreover the measure $\ms$ is invariant for
this flow and is unique. Finally, we prove that the
deterministic inviscid equation has a $\ms$-stationary solution
(for any $\nu>0$).
\end{abstract}

\section{Introduction}
The goal in this paper is to study a  class of mathematical models
related to 2D fluids. We will deal with an abstract stochastic evolution
equation in a Hilbert space of the following form
\begin{equation}\label{Abstract_Model}
du(t)+\left[\nu Au(t)+B(u(t),u(t))\right]dt=\sqrt Q dw(t),
\end{equation}
where $w$ is a cylindrical Wiener process and
$Q$ is a linear operator.
The unbounded linear operator $A$ and the bilinear operator
$B$ will satisfy certain properties related to 2D fluids that
will be given in details in the following sections.
The coefficient $\nu\ge 0$ is the viscosity.
There is an extensive literature about the existence and uniqueness of
solutions with initial data of finite energy.
Its long time behavior has also been extensively studied, including
the existence and uniqueness of invariant measures (see, e.g.,
\cite{f} and the reference therein).
In the present paper, we are
interested in the qualitative behavior of these invariant measures. In
particular, we prove the existence and uniqueness of invariant
measures of Gaussian type for the viscous case \eqref{Abstract_Model};
moreover, this Gaussian measure is proved to be invariant also
for the deterministic and inviscid  model ($\nu=0, Q=0$).

We point out that the Gaussian invariant measure that we consider here is
not that one considered in previous papers \cite{ac}, \cite{dapdeb},
\cite{deb}, \cite{af}, \cite{BF}, but has a more regular support.
In particular, the support of this measure is a Sobolev space of positive exponent.

As far as the content of this paper, in Section 2 we introduce the operators associated to the model
\eqref{Abstract_Model} with their properties
and the Gibbs measures $\ms$. We introduce the Ornstein-Uhlenbeck
equation with a suitable noise and prove that $\ms$ is its unique invariant.
In Section 3, we deal with the viscous stochastic case;
we prove the existence and uniqueness of strong solutions
and that $\ms$ is its unique invariant measure.
The uniqueness of the invariant measure is proved by means of Girsanov
theorem.
Moreover, some ergodic properties
of this measure with its rate of convergence are shown.
In Section 4, we introduce a particular example, shell
models of turbulence with an emphasis on the Sabra model.
The coefficient $\beta$ characterizing the measure $\ms$ will be related
to the coefficients $a$ and $b$ of the Sabra model through the condition
\eqref{condizbeta}. Section 5 is devoted to the deterministic inviscid
model, in particular we present our results for the inviscid Sabra model
with $\beta=1$.
For any $\nu>0$ we prove the existence of a stationary process whose law at any fixed time is $\msu$.

\section{Introduction to the model and functional setting}
\subsection{Operators and spaces}
Let $(H, |\cdot|)$ be a real separable Hilbert space
endowed with an inner product denoted by $(\cdot,\cdot)$,
and $A$ an unbounded self-adjoint positive linear operator
on $H$ with compact resolvent. We denote by
$0<\lambda_{1}\le \lambda_{2}\leq \dots$ the eigenvalues of $A$ and
by $e_{1}, e_{2}, \dots$ a complete orthonormal system in $H$
given by the eigenfunctions of the operator $A$
\[
 Ae_n=\lambda_n e_n
\]
We assume that $\lim_{n\to \infty}\lambda_n=\infty$.

Let $H_n=span\{e_1,e_2, \ldots, e_n\}$ and $\Pi_n$ the projector
operator onto $H_n$.

For any $\alpha \in \mathbb R$
we can define the power operators $A^\alpha$ as
\[
 A^\alpha x=\sum_{n=1}^\infty \lambda_n^{\alpha}(x,e_n)e_n \qquad
D(A^\alpha)=\{x=\sum_{n=1}^\infty x_n e_n:
     \sum_{n=1}^\infty \lambda_n^{2\alpha}x_n^2< \infty\}.
\]
We set
$$H^{\alpha}=D(A^{\alpha/2}).$$
Each $H^\alpha$ is a Hilbert space with scalar product
$\langle u,v\rangle_{H^\alpha}:=(A^{\alpha/2}u,A^{\alpha/2}v)$. We
denote by $\|\cdot\|_\alpha$ the norm in $H^\alpha$.

Let $B:H\times H \to H^{-1}$ be a bilinear operator; we assume that  there
exists a positive constant $c$ such that
\begin{equation}\label{Bbase}
 \|B(u,v)\|_{-1}\le c  |u| |v|.
\end{equation}

We consider the finite dimensional approximation of the bilinear
operator $B$; this is the bilinear operator $B^M$ defined as
\[
 B^{M}(u,v)=\Pi_M B(\Pi_M u, \Pi_M v)
\]
for any $M \in \mathbb N$.
For each $B^M$  we have the same estimate as \eqref{Bbase}
(with the constant $c$ independent of $M$).

For any $\nu>0$ and $\beta>0$,
let $\ms$ be the Gaussian measure $\mathcal N(0, {\frac 1 \nu}
A^{-\beta})$ (see, e.g., \cite{kuo}, \cite{dpz}).

\subsection{Assumptions}

Besides the basic properties of the operators $A$ and $B$ given above,
we present other important assumptions.
\\[1mm]
{\bf Condition (C1):} For any $\nu>0$, the operator $\nu A$
generates an analytic semigroup of contractions in $H$
and for any $p>0$ there exists $c_{p,\nu}>0$ such that
\begin{equation} \label{semigr}
 |A^p e^{-\nu A t}x|\le  \frac {c_{p,\nu}}{t^{p}} |x| \qquad \forall
 t>0, x \in H.
\end{equation}

\noindent
{\bf Condition (C2):} The bilinear operator $B$ satisfies the following properties:\\
(i) $\langle B(u,v),w\rangle=- \langle B(u,w),v\rangle$\\
(ii) $\langle B(u,v),v\rangle=0$\\
(iii) $\exists \ \beta >0$ such that $\langle B(u,u), A^{\beta}u\rangle=0$\\
for any $u,v,w$ giving meaning to the above relationships.\\

\noindent
{\bf Condition (C3):} There exists $\alpha \in [0,\beta)$
(with  $\beta$ given by {\bf (C2 iii)})
such that the embedding
$H^\beta \subset H^\alpha$ is Hilbert-Schmidt, i.e.
$$
 \sum_{n=1}^\infty \lambda_n^{\alpha-\beta}<\infty.
$$

\noindent
{\bf Condition (C4)}: for $\alpha$ and $\beta$ given in {\bf (C2)-(C3)},
$B:H^\alpha\times H^\alpha \to H^{\beta-1}$ is a continuous operator, i.e.
\begin{equation}\label{Bgenerale}
 \|B(u,v)\|_{\beta-1}\le c \|u\|_{\alpha}\|v\|_{\alpha}
 \qquad  \forall u,v \in  H^{\alpha}.
\end{equation}
Moreover, if $\alpha>0$ we assume
\begin{equation}\label{Ba}
 \|B(u,v)\|_{\alpha-1}\le c |u| \|v\|_{\alpha}
 \qquad  \forall u \in H,v \in  H^{\alpha}.
\end{equation}

\noindent
{\bf Condition (C5)}: For each $n$ set $B_n(u,v)=\langle B(u,v),e_n\rangle$.
Then we have
$$\int |B_n(x,x)|^2 \ \ms(dx)<\infty\qquad \forall n$$
and
$ B_n(x,x)$ independent of $x_n$ (where $x=\sum_n x_n e_n$).
Moreover, we require $\beta \le 1$ and 
\begin{equation}\label{condE}
 \lim_{M\to \infty} \sum_{n=1}^M \int
  |\langle B^{M}(x,x)-B(x,x), e_n\rangle|^2 \ms(dx)=0.
\end{equation}

\begin{remark}
(i) We have the relationships corresponding to assumption {\bf (C2)}:
\begin{equation}
 (B^M(u,v),w)=- (B^M(u,w),v)
\end{equation}
\begin{equation}
 (B^{M}(u,v),  v)=0
\end{equation}
\begin{equation}
 (B^{M}(u,u), A^{\beta} u)=0
\end{equation}
\\
(ii) By means of the bilinearity and of estimate \eqref{Bgenerale} we
have
\[
 \lim_{M\to \infty} \|B^M(u,v)-B(u,v)\|_{\beta-1}=0 \qquad
  \forall u, v \in H^\alpha
\]
\\
(iii) Since $\alpha \ge 0$, the inequality \eqref{Ba} implies
\begin{equation}\label{alpha-alpha}
 \|B(u,v)\|_{\alpha-1}\le c \|u\|_{\alpha}\|v\|_{\alpha}
 \qquad  \forall u,v \in  H^{\alpha}.
\end{equation}
Moreover,
\begin{equation}\label{convBM-B}
 \lim_{M\to \infty} \|B^M(u,v)-B(u,v)\|_{\alpha-1}=0 \qquad
  \forall u, v \in H^\alpha
\end{equation}
\\
(iv) Assumption {\bf (C3)} implies that the space $H^\alpha$ has full measure $\ms$,
i.e.
$\ms(H^\alpha)=1$. However, for Gaussian measures in infinite
dimensional spaces we have $\ms(H^\beta)=0$ (see, e.g., \cite{kuo}).
\end{remark}

We denote by $\mathcal L^p(\ms)$ the space of measurable functions
$\phi$ defined in the support of the measure $\ms$ and such that
$\int |\phi|^p d\ms<\infty$.

\subsection{The equations}\label{S-eqs}
Set $Q=2A^{1-\beta}$ in \eqref{Abstract_Model}, that is
we consider the following nonlinear stochastic equation
\begin{equation}\label{sequu}
du(t)+[\nu Au(t)+B(u(t),u(t))]dt = \sqrt {2A^{1-\beta}}  dw(t).
\end{equation}
In addition we deal with  the inviscid and deterministic equation
\begin{equation}\label{eul}
 \frac{du}{dt}(t)+B(u(t),u(t))=0
\end{equation}
and with the viscous linear stochastic equation
\begin{equation}\label{sequz}
 dz(t)+\nu Az(t)\ dt = \sqrt {2A^{1-\beta}} dw(t).
\end{equation}

Relationship (ii) in Assumption {\bf (C2)}
implies a formal law of conservation of energy
$E(t)=\frac 12 |u(t)|^2 $ in equation \eqref{eul}. We recall that the
energy is a conserved quantity in the motion of incompressible
inviscid fluids.

Relationship (iii) in Assumption {\bf (C2)} implies that
$S_\beta(t)=\frac 12 \|u(t)\|_\beta^2$
is a conserved quantity for equation \eqref{eul}, that is
formally we have
$$
 \dfrac {dS_\beta}{dt}(t)
 =(\dot u(t), A^{\beta} u(t))
 = - (B(u(t),u(t)), A^{\beta} u(t)) =0.
$$
For $\beta=1$, $S_1$ is the enstrophy which is a conserved quantity
in the motion of 2D incompressible inviscid fluids.

The Gaussian measure
$\ms=\mathcal N(0, {\frac 1 \nu} A^{-\beta})$ can be described
heuristically as
\[
 \ms(du)=''\frac 1Z e^{-\nu S_\beta(u)}du ''
\]
where $Z$ is a normalization constant to make $\ms$
to be a probability measure. Therefore it makes sense to see if the
measure $\ms$, described by means of the invariant quantity $S_\beta$,
is a stationary statistical solution for the inviscid equation
\eqref{eul}. To this end, we will first prove that $\ms$ is a
stationary measure for the viscous and stochastic equation
\eqref{sequu}
looking for a dynamics in the space $H^\alpha$ of full measure $\ms$.
However, the basic stochastic case to deal with is the linear equation
\eqref{sequz} for which we recall well known properties  (see \cite{dpz}).
\begin{proposition}
Let assumptions {\bf (C1)}, {\bf (C2 iii)} and {\bf (C3)}
be satisfied.\\ Then,
for any $z(0) \in H^{\alpha}$ there exists a unique strong solution
to equation \eqref{sequz} such that
\[
 z \in C([0,T];H^{\alpha})\qquad \po-a.s.
\]
\end{proposition}
The stationary process solving equation \eqref{sequz} is
\[
  \zeta(t)=\sqrt {2} \int_{-\infty}^t e^{-\nu (t-s)A } A^{\frac{1-\beta}2} dw(s)
\]
and the law of $\zeta(t)$ is $\ms$ for any time $t$.

\section{Stochastic viscous models}
We consider equation \eqref{sequu}; first we prove that there exists a
unique solution for any initial data in $H^\alpha$. The
solution is strong in the probabilistic sense and uniqueness is in
pathwise sense. Moreover, we show that $\ms$ is the unique invariant
measure associated to this stochastic equation.

\subsection{Strong solution}\label{sol-forte}
We look for dynamics in the state space $H^{\alpha}$ with $0\le
\alpha<\beta$ fulfilling assumptions {\bf (C1)}-{\bf (C4)}.
We consider any finite time interval $[0,T]$.
\begin{theorem}\label{strong_solution}
Let assumptions {\bf (C1)},{\bf (C2)}, {\bf (C3)} and
{\bf (C4)} be satisfied.\\
Then,
for any $u(0)\in H^{\alpha}$,
there exists a unique solution $u$ to equation \eqref{sequu} such
that
\[
 u \in C([0,T];H^{\alpha}) \qquad \po-a.s.
\]
Moreover, the process $u$ is a Markov process, Feller in $H^{\alpha}$.
\end{theorem}

We divide the proof in three steps in the following subsections.

\subsubsection{Existence of strong solutions}
We use a well known trick to study a stochastic semilinear equation
with additive noise:  we set $v=u-z$. Then
\begin{equation}\label{sequv}
 \frac {dv}{dt}(t)+\nu Av(t)+ B(v(t)+z(t),v(t)+z(t))=0
\end{equation}
with $v(0)=u(0)-z(0)$. Set $z(0)=0$.

\begin{proposition}
We consider the same assumptions as in Theorem
\ref{strong_solution}.
Let $v(0)\in H^{\alpha}$.
Then there exists a solution to equation \eqref{sequv} such
that
\[
 v \in C([0,T];H^{\alpha})\cap L^2(0,T;H^{1+\alpha}) \qquad \po-a.s.
\]
\end{proposition}
\proof
We proceed pathwise. Take the scalar product of the left hand side of
equation \eqref{sequv} with $v$ in $H$; we get
some a priori estimates
\[\begin{split}
\frac 12 \frac {d}{dt} |v|^2+ \nu
\|v\|^2_1
&=- \langle B(v+z,v+z),v\rangle
\\&=
  - \langle B(v+z,z),v\rangle \quad\text{ by  {\bf (C2 ii)}}
\\& \le
  \|B(v+z,z)\|_{-1}\|v\|_{1}
\\& \le
  c|v+z| |z| \|v\|_{1} \quad\text{ by }\eqref{Bbase}
\\& \le
  \frac \nu 2 \|v\|_{1}^2+ \frac{c_\nu}2 |z|^2 |v|^2
   + \frac{c_\nu}2 |z|^4
\end{split}
\]
by Young inequality,
for some positive constant $c_\nu$.
Henceforth, we denote by $c_\nu$ a generic constant depending on
$\nu$.

Therefore
\begin{equation}\label{primastima}
 \frac {d}{dt} |v|^2+ \nu \|v\|^2_1
 \le
 c_\nu |z|^2 |v|^2  +c_\nu |z|^4.
\end{equation}
Hence,  Gronwall inequality applied to
\[
 \frac {d}{dt} |v|^2 \le  c_\nu |z|^2 |v|^2  +c_\nu |z|^4
\]
gives
\begin{equation}\label{stime0}
 \sup_{0\le t \le T}|v(t)|^2\le
  e^{c_\nu T\|z\|_{C([0,T];H)}^2}
 \left( |v(0)|^2 +  c_\nu T \|z\|_{C([0,T];H)}^4\right)
  <\infty
\end{equation}
and integrating in time \eqref{primastima}
\begin{equation}\label{stime01}
 \nu \int_0^T  \|v(s)\|_{1}^2 ds\le
  |v(0)|^2+ T c_\nu
  \left( \|z\|_{C([0,T];H)}^2 \|v\|_{C([0,T];H)}^2+ \|z\|_{C([0,T];H)}^4 \right)
 <\infty.
\end{equation}
Moreover, when $\alpha \geq 0$ we proceed in a similar way:
we take the scalar product of the left hand side of equation
\eqref{sequv} with $A^\alpha v$ in $H$; then
\[\begin{split}
\frac 12 \frac {d}{dt} \|v\|^2_{\alpha}+ \nu \|v\|^2_{1+\alpha}
&=-\big(A^{\frac{-1+\alpha}2} B(v+z,v+z),A^{\frac{1+\alpha}2} v\big)
\\& \le
  \|B(v+z,v+z)\|_{-1+\alpha}\|v\|_{1+\alpha}
\\& \le
  c|v+z| \|v+z\|_{\alpha} \|v\|_{1+\alpha}
            \quad\text{ by }\eqref{Ba}
\\& \le
 \frac \nu 2 \|v\|_{1+\alpha}^2
 + \frac {c_\nu}2 (|v|^2 + |z|^2) \|v\|^2_{\alpha}
 + \frac {c_\nu}2 \|z\|^4_{\alpha}.
\end{split}
\]
This gives
\begin{equation}\label{stime3}
 \frac {d}{dt} \|v\|^2_{\alpha}+ \nu \|v\|^2_{1+\alpha}
\le
 c_\nu \Big( |v|^2 + \|z\|_\alpha^2\Big)  \|v\|^2_{\alpha}+ c_\nu \|z\|^4_{\alpha}.
\end{equation}
Therefore, using \eqref{stime0} and the fact that $\alpha\geq 0$ we get
\begin{multline}\label{L-inf-alpha}
 \sup_{0\le t \le T} \|v(t)\|^2_{\alpha}
 \le \|v(0)\|^2_{\alpha}  e^{c_\nu\int_0^T  (|v(t)|^2 + \|z(t)\|^2_{\alpha})dt}
\\    + c_\nu \int_0^T e^{c_\nu \int_t^T (|v(s)|^2 + \|z(s)\|^2_{\alpha})ds}
       \|z(t)\|^4_{\alpha} dt
 < \infty
\end{multline}
and integrating in time \eqref{stime3}
\[
 \int_0^T  \|v(s)\|^2_{1+\alpha} ds < \infty.
\]

Actually, the a priori estimates are for the Galerkin approximation
$v^M$. We define the  Galerkin problem associated to \eqref{sequu}
\begin{equation}\label{gal}
\begin{cases}
du^M(t)+[\nu Au^M(t)+B^{M}(u^M(t),u^M(t))]dt
      = \Pi_M \sqrt {2A^{1-\beta}}\ dw(t)\\
 u^M(0)=\Pi_M x
\end{cases}
\end{equation}
where $M$ is any positive integer.
Similarly we have
\[
 \frac{dv^M}{dt}(t)+ \nu Av^M(t)+ B^M(v^M(t)+z^M(t),v^M(t)+z^M(t))=0
\]
with $z^M(t)=\Pi_M z(t)$.

The previous estimates give
\begin{equation}\label{Minf}
 \sup_M \|v^M\|^2_{L^\infty(0,T;H^{\alpha})} < \infty
\end{equation}
\begin{equation}\label{M2}
 \sup_M \|v^M\|^2_{L^2(0,T;H^{1+\alpha})} < \infty
\end{equation}
In addition $\dfrac{dv^M}{dt}$ is bounded: indeed
\[
 \frac{dv^M}{dt}(t) = - \nu Av^M(t)- B(v^M(t)+z^M(t),v^M(t)+z^M(t));
\]
using \eqref{M2}-\eqref{Minf}, we have that
the first term in the r.h.s belongs to the space $ L^2(0,T;H^{\alpha-1})$
and the second to the space $C([0,T];H^{\alpha-1})$
(use \eqref{alpha-alpha}) and thus in $L^2(0,T;H^{\alpha-1})$.
Then
\begin{equation}\label{Mder}
  \sup_M \|\frac{dv^M}{dt}\|^2_{L^2(0,T;H^{\alpha-1})} < \infty.
\end{equation}
Since the space $\{v: v \in L^2(0,T;H^{1+\alpha}),\frac{dv}{dt}
 \in L^2(0,T;H^{\alpha-1})\}$
is compactly embedded in the space $L^2(0,T;H^{\alpha})$,
from \eqref{Minf}-\eqref{Mder}
we get that there exists a subsequence $\{v^{M_i}\}$ weakly convergent to a $v$
in $L^2(0,T;H^{1+\alpha})$, weakly-* convergent in
$L^\infty(0,T;H^{\alpha})$ and
strongly convergent in $L^2(0,T;H^{\alpha})$.
By means of the bilinearity of $B$, of the strong convergence result
and of \eqref{convBM-B}, we
conclude that the limit $v$ fulfils \eqref{sequv}.

The fact that $v \in C([0,T];H^{\alpha})$
comes from a result in Temam \cite{Temam} (Lemma 1.4. page 263): if
$v\in L^2(0,T;H^{1+\alpha})$ and $\frac{dv}{dt} \in L^2(0,T;H^{-1+\alpha})$,
then $v \in C([0,T];H^{\alpha})$.
\hfill $\Box$

\begin{remark}
We can prove also the uniqueness of this solution $v$, but we do not
need it here. Anyway, the proof of uniqueness would be based on the
same estimes as in the next Section \ref{Suni}.
\end{remark}

\medskip
We conclude for $u=v+z$.
\begin{proposition}
We consider the same assumptions as in Theorem
\ref{strong_solution}.
Let $u(0)\in H^{\alpha}$.
Then there exists a solution to equation \eqref{sequu} such
that
\[
 u \in C([0,T];H^{\alpha}) \qquad \po-a.s.
\]
\end{proposition}

\subsubsection{Pathwise uniqueness}\label{Suni}
Now we prove that the strong solution $u$ constructed in the previous
section is pathwise unique, that is
\begin{proposition}
We consider the same assumptions as in Theorem
\ref{strong_solution}.
Let  $u_1, u_2$ be two solutions to equation \eqref{sequu} with the
same initial data,
defined on the same stochastic basis and with the same Wiener process.
Then
$u_1=u_2$ $\po$-a.s., the equality being in $C([0,T];H^{\alpha})$.
\end{proposition}
\proof
We proceed pathwise.
Let $u_1, u_2 \in C([0,T]; H^{\alpha})$ be two paths (for fixed
$\omega$ in a set of $\po$-measure 1).
\\
Set $U=u_1-u_2$.
Then $U\in C([0,T]; H^{\alpha})$ and it solves an equation which is
deterministic (for any path):
\begin{equation}\label{U1}
 \frac{dU}{dt}+\nu AU + B(u_1,u_1)-B(u_2,u_2)=0; \qquad U(0)=0.
\end{equation}
First, we notice that $U$ is more regular than the $u_i$'s (the noise
term has desappeared and we expect more regularity as
for equation \eqref{sequv}).

By the bilinearity of the operator $B$, we have
\begin{equation}\label{U2}
 \frac{dU}{dt}+\nu AU + B(u_1,U)+B(U,u_2)=0; \qquad U(0)=0.
\end{equation}
We get an a priori estimate:
\[
\begin{split}
 \frac 12 \frac{d}{dt}\|U(t)\|^2_{\alpha} &+ \nu \|U(t)\|^2_{1+\alpha}
\\
 &=
 -\left( A^{\frac{\alpha-1}2}[B(u_1(t),U(t))+B(U(t),u_2(t))],
 A^{\frac{\alpha+1}2}U(t)\right)
\\
&\le
 \big[\|u_1(t)\|_{\alpha}+\|u_2(t)\|_{\alpha}\big]
  \|U(t)\|_{\alpha}\|U(t)\|_{1+\alpha}
\text{ by } \eqref{Bgenerale}
\\
&\le
\frac \nu 2 \|U(t)\|_{1+\alpha}^2+\frac{c_\nu}2
 \big[\|u_1(t)\|^2_{\alpha}+\|u_2(t)\|^2_{\alpha}\big]
  \|U(t)\|^2_{\alpha}.
\end{split}
\]
Therefore
\[
  \frac{d}{dt}\|U(t)\|^2_{\alpha}\le
  c_\nu  \big[\|u_1(t)\|^2_{\alpha}+\|u_2(t)\|^2_{\alpha}\big]
  \|U(t)\|^2_{\alpha};
\]
from this, by Gronwall inequality  follows
\begin{equation}\label{stima-differ}
 \|U(t)\|^2_{\alpha}\le \|U(0)\|^2_{\alpha}
 e^{c_\nu\int_0^t [\|u_1(s)\|^2_{\alpha}+\|u_2(s)\|^2_{\alpha}]ds}.
\end{equation}
Finally, $U(t)=0$ for all $t$, since $U(0)=0$.

\begin{remark}
Markovianity is inherited from the Galerkin approximations.
\end{remark}

\subsubsection{Feller property}
Let us denote by $u(t;x)$ the solution of equation \eqref{sequu} with
initial data $x$, by $B_b(H^{\alpha})$  the space of Borel bounded functions
$\phi:H^{\alpha} \to \mathbb R$ and by $C_b(H^{\alpha})$ its  space of 
continuous bounded functions.

Define the Markov semigroup $P_t:B_b(H^{\alpha})\to B_b(H^{\alpha})$ as
\[
 P_t\phi(x)=\mathbb E[\phi(u(t;x))].
\]

This is a contraction semigroup.
Moreover, it is Feller in $H^{\alpha}$, that is
\[
 P_t:C_b(H^{\alpha})\to C_b(H^{\alpha}).
\]
This comes from the estimates for the pathwise
uniqueness.
Indeed, if $\|x-y\|_{\alpha}\to 0$ then \eqref{stima-differ} gives
\begin{equation}\label{stimaF}
  \|u(t;x)-u(t;y)\|^2_{\alpha}\le \|x-y\|^2_{\alpha}
 e^{c_\nu\int_0^t [\|u(s;x)\|^2_{\alpha}+\|u(s;y)\|^2_{\alpha}]ds}
\end{equation}
for $t>0$ fixed.
By \eqref{L-inf-alpha} we get a uniform estimate of
$\|u(\cdot;x)\|^2_{L^\infty(0,T;H^{\alpha})}$ when $\|x\|_{\alpha}$ is
bounded, i.e.
\[
 \forall R>0\; \exists C_R:
\sup_{\|x\|_{\alpha}\le R}\|u(\cdot;x)\|_{C([0,T];H^\alpha)}<C_R.
\]
Hence, when $\|x-y\|_\alpha \to 0$ from \eqref{stimaF} we get
$\|u(t;x)-u(t;y)\|_{\alpha}\to 0$. We conclude that
$\phi(u(t;x)) \to \phi(u(t;y))$ for $\phi\in C_b(H^{\alpha})$ and therefore
$\mathbb E[\phi(u(t;x))] \to \mathbb E[\phi(u(t;y))]$
by the dominated convergence.
This means that $P_t\phi \in C_b(H^\alpha)$ for any $t>0$ and
$\phi\in C_b(H^\alpha)$.

\subsection{Invariant measure}\label{Sminv}
We prove the following theorem:

\begin{theorem}\label{Uniqueness_measure}
Besides the assumptions of Theorem \ref{strong_solution}
we consider {\bf (C5)}.
Then, $\ms$ is the unique invariant measure for equation \eqref{sequu},
that is
\begin{equation}\label{invarianza}
 \int P_t \phi\ d\ms= \int \phi \ d\ms \qquad
\forall \phi  \in \mathcal L^1(\ms) \text{ and } t\ge 0.
\end{equation}
\end{theorem}

First, we show that $\ms$ is an invariant measure
for the nonlinear equation \eqref{sequu}. Then we prove
that this is indeed the unique invariant measure.

A consequence of this result is the following
\begin{corollary}\label{e-staz}
Given any initial data with law $\ms$, there exists a unique
stationary solution of equation \eqref{sequu}
whose law at any fixed time is $\ms$.
\end{corollary}

Finally, in Section \ref{la-rate}
we analyse the rate of convergence of $P_t \phi$, as $t \to \infty$. 

To prove our results, we need to introduce the Kolmogorov operator
associated to the stochastic equation \eqref{sequu}. Let
$FC^\infty_b$ be  the space of infinitely differentiable cylindrical
functions  bounded and with bounded derivatives;
$\phi \in FC^\infty_b$
means that there exist $ m \in \mathbb N$, $\tilde \phi \in
C^\infty_b(\mathbb R^m)$ and multiindices $(i_1, i_2, \ldots, i_m)$
such that
\[
 \phi(x)=\tilde \phi\big((x,e_{i_1}), (x,e_{i_2}), \ldots, (x,e_{i_m})\big).
\]
We set $\dfrac{\partial \phi}{\partial x_i}
= \dfrac{\partial \tilde \phi}{\partial x_i}$ with $x_i=(x,e_{i})$.
 $FC^\infty_b$ is a dense
subset of  $\mathcal L^p(\ms)$ for any $p\ge 1$.\\
We define the Kolmogorov operator first on these very regular
functions $\phi \in FC^\infty_b$ as
\begin{equation}\label{Kinv}
K\phi(x)=  \sum_j\Big[
  \lambda_j^{1-\beta}\frac{\partial^2 \phi}{\partial x_{j}^2}(x)
    - B_{j}(x,x) \frac{\partial \phi}{\partial x_{j}}(x)
    - \nu \lambda_j x_{j}\frac {\partial \phi}{\partial x_{j}}(x)\Big].
\end{equation}
We have that  $K\phi\in \mathcal L^1(\ms)$ for any $\phi \in  FC^\infty_b$
(use that each $B_j \in \mathcal L^1(\ms)$ and the sums are finite).

\subsubsection{Existence of the invariant measure}\label{321}
We know that the linear stochastic equation
\eqref{sequz} has $\ms$ as unique invariant measure, that is $\ms$ is
the unique probability measure such that
\[
 \int \mathbb E [\phi(z(t;x))]\ms(dx)= \int \phi(x) \ms(dx)
 \qquad \forall t\ge 0, \phi \in B_b(H^\alpha)
\]
(see \cite{dpz,dpz2}).
Actually we can define the latter relationship for all
$\phi \in \mathcal L^p(\ms)$, given any $1\le p<\infty$
(see, e.g., \cite{daprato_birkhauser,daprato_springer}).

Now, we want to show that $\ms$ is an invariant measure also
for the nonlinear equation \eqref{sequu}.
 The role of the nonlinear term $B$ is analyzed first
considering the finite dimensional $B^{M}$ and then
passing to the limit as $M\to \infty$. Here we need
\eqref{condE} of {\bf (C5)}.

First, we prove that $\ms$ is an {\it infinitesimally}
invariant measure for equation \eqref{sequu}
in the sense that
\begin{equation}\label{KM-inv}
 \int K\phi\ d\ms = 0 \quad \forall \phi \in  FC^\infty_b .
\end{equation}
Indeed, we can write $K$ as the sum of two operators,
$K=Q+L$, with domains $FC^\infty_b$
and we have the infinitesimal invariance for both these operators.
We integrate by parts:
\begin{equation}\label{invQ}
\int Q \phi\  d\ms \equiv
\int  \sum_{j} \Big[\lambda_j^{1-\beta}\frac{\partial^2 \phi}{\partial
    x_{j}^2}(x)-
 \nu \lambda_j x_{j}\frac {\partial \phi}{\partial
      x_{j}}(x)\Big]\ \ms(dx)=0
\end{equation}
and
\begin{multline}\label{invL}
\int L \phi\  d\ms \equiv
  - \int  \sum_{j} B_{j}(x,x) \frac{\partial \phi}{\partial x_{j}}(x)\ \ms(dx)
\\
= -\nu \int
 \underbrace{\sum_{j}\lambda_j^{\beta}
             B_{j}(x,x) x_{j}}_{=0 \;\;\text{ by } {\bf (C2 iii)}}
 \ \phi(x)\ \ms(dx)=0
\end{multline}
since $B_j$ does not depend on the variable $x_j$.

With similar computations,
we get that the Kolmogorov operator $(K,FC^\infty_b)$ is
dissipative in $\mathcal L^1(\ms)$, that is
\[
  \int \phi \ K\phi \ d\ms\le 0 \qquad \forall \phi \in  FC^\infty_b
\]
(see also \eqref{Kdissip}). Hence it is closable (see \cite{pazy}).

Now we use an approximative criterium of Eberle \cite{e}
in order  to show that the measure $\ms$ is an
invariant measure for equation \eqref{sequu}:
\begin{proposition}\label{la-dimo-E}
Besides the assumptions of Theorem \ref{strong_solution}
we consider {\bf (C5)}.
Then, 
the closure operator $\overline K$ 
of the Kolmogorov operator $(K,FC^\infty_b)$ in $\mathcal L^1(\ms)$
generates a sub-Markovian strongly continuous semigroup
$T_t=e^{\overline K t}$  in $\mathcal L^1(\ms)$.
\end{proposition}
Moreover, $T_t$ is the only strongly continuous semigroup on $\mathcal L^1(\ms)$
which has generator that extends $(K, FC^\infty_b)$ (see Appendix A in
\cite{e}).

We postpone the proof of this result and continue our analysis.
Since $FC^\infty_b$ is a core for the infinitesimal
generator of $T_t$  in $\mathcal L^1(\ms)$
by density  from \eqref{KM-inv} we get that
\[
 \int \overline K\phi\ d\ms = 0 \qquad \forall \phi \in D(\overline K).
\]
This is equivalent to
\begin{equation}\label{Tinv}
 \int T_t \phi\ d\ms= \int \phi \ d\ms \qquad
\forall \phi  \in \mathcal L^1(\ms) \text{ and } t\ge 0.
\end{equation}

Now, we go back to the semigroup $\{P_t\}$;
it has been constructed by means of the unique solution $u$ of
equation \eqref{sequu} such that
$u(t;x) \in H^\alpha$
(for any $t>0$, $x \in H^\alpha$). On the other hand, the analytical analysis
of the Kolmogorov operator has led to the construction of the
semigroup $\{T_t\}$; it provides a martingale solution to the
stochastic equation \eqref{sequu}
(see, e.g., \cite{e} and references therein).
By our previous results of Section \ref{sol-forte} on the stochastic
equation \eqref{sequu} we can relate these semigroups
and get that  the semigroup $\{P_t\}$ can be extended to
$\mathcal L^1(\ms)$ (where this semigroup is exactly
$\{T_t\}$).

 Henceforth, we denote these semigroups in $C_b(H^\alpha)$ and
$\mathcal L^1(\ms)$ with the same symbol $P_t$. Therefore \eqref{Tinv}
completes our proof of \eqref{invarianza}.

Now, we go back to the proof of Proposition \ref{la-dimo-E}.
We refer to \cite{e} for all the details; 
in particular, we use Theorem 5.2, Corollary 5.3, Lemma 5.11
and (5.46) at page 226 of \cite{e} with $p=1$.
\\
\proof (of Proposition \ref{la-dimo-E})
From Lumer-Phillips theorem we know that the closure of the operator 
$(K,FC^\infty_b)$ in $\mathcal L^1(\ms)$ generates a strongly
continuous semigroup $T_t$ if and only if the range of 
$(\lambda-K, FC^\infty_b)$
is dense in $\mathcal L^1(\ms)$ for some (and all) $\lambda>0$.
To prove the density result, we use an approximative criterium:
\begin{equation}\label{approx-cr}
\forall F \in \mathcal L^1(\ms) \;\forall \varepsilon>0\qquad
\exists \ v \in FC^\infty_b:\;
\|(\lambda-K)v-F\|_{\mathcal L^1(\ms)}<\varepsilon.
\end{equation}
Now, we take  $F\in \mathcal L^1(\ms)$. Then there
exists a sequence $\{F^N\}_{N \in \mathbb N}$ 
with $F^N \in C^\infty_b(\mathbb R^N)$ and
\begin{equation}
 \lim_{N \to \infty} \|F^N-F\|_{\mathcal L^1(\ms)}=0, \qquad
 \sup_N \|F^N\|_{C_b}<\infty.
\end{equation}
On the other hand, the assumption 
$B_n \in \mathcal L^2(\ms)$ (for any $n$) implies 
that $B^N \in \mathcal L^2(\ms)$ for any $N$, and therefore
there exists a sequence $\{C^N\}_{N \in \mathbb N}$ with 
$C^N \in C^\infty_b(\mathbb R^N \to \mathbb R^N)$ and
\begin{equation}
 \|B^N-C^N\|_{\mathcal L^2(\ms)}\le \frac 1N.
\end{equation}
Bearing in mind \eqref{condE}, this implies that 
\begin{equation}\label{BN-CN}
 \|\Pi_N B-C^N\|_{\mathcal L^2(\ms)}\le 
 \|\Pi_N B-B^N\|_{\mathcal L^2(\ms)}+\|B^N-C^N\|_{\mathcal L^2(\ms)} \to 0 
\end{equation}
as $N \to \infty$.

For each $N$, we introduce a regularized finite dimensional
 Kolmogorov operator $K^N$ acting on
functions $\phi \in C_b^\infty(\mathbb R^N)$:
\[
 (K^N\phi)(x)= \sum_{j=1}^N\Big[
  \lambda_j^{1-\beta}\frac{\partial^2 \phi}{\partial x_{j}^2}(x)
    - C^N_{j}(x) \frac{\partial \phi}{\partial x_{j}}(x)
    - \nu \lambda_j x_{j}\frac {\partial \phi}{\partial x_{j}}(x)\Big].
\]
It as smooth coefficients.
Therefore, given $F^N \in C^\infty_b(\mathbb R^N)$ and $\lambda>0$ the equation
\begin{equation}\label{autoval}
 (\lambda-K^N)\phi^N=F^N
\end{equation}
has a unique solution $\phi^N \in C^\infty_b(\mathbb R^N)$; moreover
\begin{equation}\label{propSubM}
 \lambda\|\phi^N\|_{C_b} \le \|F^N\|_{C_b} .
\end{equation}
Further, setting $|\sqrt{A^\gamma} D\phi^N|^2=
\sum_{j=1}^N \lambda_j^\gamma |\frac{\partial \phi^N}{\partial x_j}|^2
$, by a straightforward computation we have
\begin{equation*}
K (\phi^N)^{2}=  2\phi^N\ K\phi^N +2|\sqrt{A^{1-\beta}} D \phi^N|^{2}.
\end{equation*}
Using the infinitesimal invariance \eqref{KM-inv}, we have
\[\begin{split}
 \|\sqrt {A^{1-\beta}} D\phi^N\|^2_{\mathcal L^2(\ms)} 
& = -\int \phi^N \ K\phi^N d\ms 
\\
& =   -\int \phi^N  [K\phi^N-K^N\phi^N] d\ms - \int \phi^N \ K^N
       \phi^N d\ms
\\
& =  -\int \phi^N \sum_{j=1}^N (B_j-C^N_j)\frac {\partial
  \phi^N}{\partial x_{j}}  d\ms
 + \int \phi^N (F^N-\lambda \phi^N) d\ms
\\
& \le
 \|\phi^N\|_{C_b} \|\Pi_N B - C^N\|_{\mathcal L^2(\ms)} \|D\phi^N\|_{\mathcal L^2(\ms)}
\\
&\qquad + \|\phi^N\|_{C_b} \big(\|F^N\|_{C_b}+\lambda \|\phi^N\|_{C_b}\big).
\end{split}\]
Using \eqref{propSubM} and the fact  that 
$\lambda_1^{1-\beta} \|D\phi^N\|^2_{\mathcal L^2(\ms)}
\le  \|\sqrt{A^{1-\beta}} D\phi^N\|^2_{\mathcal L^2(\ms)}$
for $\beta\le 1$, we get that there exists  a constant $C_{\beta,\lambda}>0$ such
that
\begin{equation}\label{stima-gradiente}
 \|D\phi^N\|_{\mathcal L^2(\ms)} \le C_{\beta,\lambda} \|F^N\|_{C_b}
  \left[ 1+  \| \Pi_N B - C^N\|_{\mathcal L^2(\ms)}\right].
\end{equation}

Let us go back to \eqref{approx-cr}; by \eqref{autoval} we have
\[
 (\lambda-K)\phi^N-F
          = (K^N-K)\phi^N+F^N-F
        \equiv \sum_{j=1}^N [B_j-C^N_j]\frac{\partial \phi^N}{\partial x_j}+F^N-F.
\]
Integrating with respect to the measure $\ms$ we get
\[\begin{split}
 \|(\lambda-K)\phi^N-F\|_{\mathcal L^1(\ms)}
& \le
  \int \Big| \sum_{j=1}^N [B_j-C^N_j]
       \frac{\partial \phi^N}{\partial x_j}\Big|d\ms
   +\|F^N-F\|_{\mathcal L^1(\ms)}\\
& \le
  \|\Pi_N B-C^N\|_{\mathcal L^2(\ms)} \|D\phi^N\|_{\mathcal L^2(\ms)} 
   +\|F^N-F\|_{\mathcal L^1(\ms)}
\end{split}
\]
by Schwarz inequality.
Using \eqref{stima-gradiente}  we find
\begin{multline*}
 \|(\lambda-K)\phi^N-F\|_{\mathcal L^1(\ms)}\le C_{\beta,\lambda}
  \|\Pi_N B-C^N\|_{\mathcal L^2(\ms)}\|F^N\|_{C_b}
   \Big[1+\|\Pi_NB-C^N\|_{\mathcal L^2(\ms)}\Big]
\\
       +\|F^N-F\|_{\mathcal L^1(\ms)}.
\end{multline*}
Bearing in mind the assumptions on the approximating terms and 
\eqref{BN-CN}, we  find \eqref{approx-cr}.
\hfill $\Box$

\begin{remark}\label{oss-Lp}
Because of the invariance of the measure $\ms$, the
contraction semigroup $P_t$ in $C_b(H^\alpha)$ can be uniquely extended to a
strongly continuous
contraction semigroup in $\mathcal L^p(\ms)$ also for any $p>1$.
Indeed,
\[
 |P_t \phi(x)|^p= |\mathbb E [\phi(u(t;x))]|^p \le
                    \mathbb E [|\phi(u(t;x))|^p] = P_t |\phi|^p(x)
\]
and by the invariance of the measure $\ms$
\[
 \int |P_t\phi|^p d\ms\le \int P_t |\phi|^p d\ms = \int |\phi|^p d\ms.
\]
Since $C_b(H^\alpha)$ is dense in $\mathcal L^p(\ms)$, we can uniquely
define the semigroup on $\mathcal L^p(\ms)$ for any $p> 1$.
We use the same symbol $P_t$ to denote all these semigroups.

Notice that in condition {\bf (C5)} we require
$\int | B_n(x,x)|^2 \ \ms(dx)<\infty$
for any $n$.
Therefore  $K:FC^\infty_b \to \mathcal L^2(\ms)$.
Moreover, according to Corollary 5.3 of \cite{e}, we have that the restriction
of $T_t$ to $\mathcal L^2(\ms)$ is a strongly continuous semigroup on
$\mathcal L^2(\ms)$ and the
generator of this semigroup again extends $(K, FC^\infty_b)$. In the
sequel we will use the same symbol to denote these semigroups in
both spaces $\mathcal L^1(\ms)$ and $\mathcal L^2(\ms)$.
\end{remark}

\subsubsection{Uniqueness of the invariant measure}
Now we prove that equation \eqref{sequu} has at most one invariant measure.
We use the results of Section \ref{sol-forte}.

Let $ R(t,x,\cdot) $ be the law of $z(t;x)$ and
$ P(t,x,\cdot) $ be the law of $u(t;x)$.
Then any $ R(t,x,\cdot) $ is equivalent to the Gibbs measure $\ms$
(see, e.g., \cite{dpz2}); we write it as $R(t,x,\cdot) \sim \ms$.
Moreover we have that
\begin{equation}\label{gir}
 \int_0^T |\sqrt{A^{\beta-1}}B(z(t),z(t))|^2 dt <\infty \qquad \po-a.s.
\end{equation}
and
\begin{equation}\label{gir2}
 \int_0^T |\sqrt{A^{\beta-1}} B(u(t),u(t))|^2 dt <\infty \qquad \po-a.s.
\end{equation}
For this use that
$\|B(x,x)\|_{\beta-1}\le c \|x\|^2_{\alpha}$ from assumption
\eqref{Bgenerale} and that $\po\{z \in C([0,T];H^\alpha)\}=
\po\{u \in C([0,T];H^\alpha)\}=1$.

According to Theorem 9.2  in \cite{Fe},
\eqref{gir}-\eqref{gir2}
imply that the measure $P(t,x,\cdot)$ is equivalent to
$R(t,x,\cdot)$.
On the other side $R(t,x,\cdot)\sim R(s,y,\cdot)\sim \ms$, hence
we get
\[
 P(t,x,\cdot) \sim  P(t,y,\cdot)\sim \ms
\]
for any $x,y \in H^{\alpha}$ and $t>0$.
Using Doob theorem (see, e.g., Theorem 4.2.1 in \cite{dpz2}), we
deduce that there exists at most one invariant measure.

By means of the existence result of the previous section,
we get that $\ms $ is the unique invariant measure
for equation \eqref{sequu}. Moreover, it is strongly mixing
\begin{equation}\label{strongly_mixing}
 \lim_{t\rightarrow\infty}P(t,x,\Gamma)=\ms(\Gamma)
\end{equation}
for arbitrary $x\in H^{\alpha}$ and  Borel set
$\Gamma$ in  $H^{\alpha}$.

\subsubsection{Rate of convergence}\label{la-rate}
Now, we consider the semigroup $P_t$ in $\mathcal L^2(\ms)$
(see Remark \ref{oss-Lp}).

We recall the ''Carr\'e du champ'' identity. For the reader's
convenience we give the proof (see, e.g.,  \cite{daprato_birkhauser})
\begin{proposition}
Besides the assumptions of Theorem \ref{strong_solution}
we consider {\bf (C5)}. Then, we have
\begin{equation}\label{cdc}
 \int \phi \ \overline K \phi\ d\ms= - \int |\sqrt
      {A^{1-\beta}}D\phi|^2d\ms \qquad
       \forall \phi \in D(\overline K).
\end{equation}
\end{proposition}
\proof
First we take $\phi \in FC^\infty_b$.
A straightforward computation yields that
\begin{equation*}
K \phi^{2}=  2\phi\ K\phi +2|\sqrt{A^{1-\beta}} D \phi|^{2}.
\end{equation*}
By the $\ms$-infinitesimal invariance, we have $\int K \phi^{2} \ d\ms=0 $; thus
\begin{equation}\label{Kdissip}
\int \phi\ K\phi \ d\ms=
-  \int |\sqrt{A^{1-\beta}} D \phi |^{2} d\ms.
\end{equation}
Now, taking $\phi\in D(\overline K)$, we use that $FC^\infty_b$ is a core for
$\overline K$; therefore there exists a sequence $\{\phi_n\}\subset
FC^\infty_b$ such that
\[
 \phi_n \to \phi, \; K \phi_n \to \overline K \phi \qquad \text{ in }
 \mathcal L^2(\ms).
\]
From \eqref{Kdissip} we get
\[
 \int |\sqrt{A^{1-\beta}} D (\phi_n-\phi_m) |^{2} d\ms\le
 \int |\phi_n-\phi_m| |K(\phi_n-\phi_m)| d\ms.
\]
Hence, the sequence $\{\sqrt{A^{1-\beta}} D \phi_n\}$
is a Cauchy sequence in $\mathcal L^2(\ms)$
and we get \eqref{cdc}.
\hfill $\Box$

Now, given  $\phi\in \mathcal{L}^{2}(\ms)$
we set $\overline \phi=\int \phi \ d\ms$; then we have the following theorem
on the rate of convergence of $P_t \phi$ as $t \to \infty$.
\begin{theorem}\label{rate_CV}
Besides the assumptions of Theorem \ref{strong_solution}
we consider {\bf (C5)}. Then
\[
\int \mathbb |P_{t}\phi(x)-\overline \phi|^2 \ms(dx)
\le
 e^{-\lambda_{1} t} \int \mathbb  |\phi(x)-\overline \phi|^2 \ms(dx)
\]
for any $\phi\in \mathcal{L}^{2}(\ms)$ and $t>0$.
\end{theorem}
\proof

Let us define the space
\begin{equation}\label{L_{0}}
\mathcal L^2_{0}(\ms)=\left\{\phi\in \mathcal L^2(\ms):\ \ \overline \phi=0\right\};
\end{equation}
it is not difficult to prove that it is invariant for the semigroup
$P_{t}$ (see \cite{daprato_springer}).

First,
let us take $\phi\in \mathcal L^2_{0}(\ms)\cap D(\overline K)$; then
$P_t \phi\in \mathcal L^2_{0}(\ms)\cap D(\overline K)$ and
by the Hille-Yosida theorem
\[
 \frac{d}{dt} P_t \phi = \overline KP_t \phi.
\]
Therefore, bearing in mind \eqref{cdc}
\[
 \frac 12 \frac {d}{dt} \int |P_t \phi|^2 \ d\ms=
 \int P_t \phi \ \overline K P_t \phi \ d\ms=
 - \int |\sqrt{A^{1-\beta}}D_x P_t \phi|^2 d\ms
\]
Since a Gaussian measure fulfils the spectral gap inequality
(see \cite{goldys})
we have
\[
  \quad \int |\sqrt{A^{1-\beta}} D_xP_t\phi(x)|^2 \ms(dx)
  \ge \frac{\lambda_1}2
  \int [P_t\phi(x)]^2 \ms(dx)
\]
where $\lambda_1>0$ is the first eigenvalue of the operator $A$.
By the two latter relationships we get
\[
 \frac {d}{dt} \int |P_t \phi|^2 d\ms\le
  -\lambda_1 \int |P_t \phi|^2 d\ms.
\]
Hence, using Gronwall lemma, we have that for any $t>0$
\begin{equation}\label{rate0}
 \int \mathbb |P_{t}\phi|^2 \ d\ms
 \leq
 e^{-\lambda_{1} t} \int \mathbb  |\phi|^2 \  d\ms
\qquad \forall \phi\in \mathcal{L}_{0}^{2}(\ms)\cap D(\overline K).
\end{equation}
Now we take $\phi \in D(\overline K)$; replacing $\phi$ with
$\phi-\overline{\phi}$ in \eqref{rate0}, we obtain that
\[
 \int \mathbb |P_{t}\phi-\overline \phi|^2 \ d\ms
 = \int \mathbb |P_{t}(\phi-\overline \phi)|^2 \ d\ms
 \le e^{-\lambda_{1} t} \int \mathbb  |\phi-\overline \phi|^2 \ d\ms.
\]
Using that $D(\overline K)$ is dense in $\mathcal L^2(\ms)$ we get the result.

\hfill $\Box$

\section{An example: shell models of turbulence}
Shell models of turbulence describe the evolution of
complex Fourier-like components of a scalar velocity field. Here we
present the details for the SABRA shell model (see \cite{sabra}),
but the same results hold for the GOY shell model (see \cite{G,goy}).
In recent years there has been an increasing interest in these
fluid dynamical models, both for the deterministic and the stochastic
case (see also \cite{clt}, \cite{bbbf},  \cite{bm},  \cite{cm}). They are
easier to analyze than the Navier-Stokes or Euler equations, but they
retain many important features of the true hydrodynamical models.

Instead of dealing with complex valued unknowns we deal with the real
and imaginary part of each component  of the scalar velocity
field (for the basic settings we follow \cite{BF});
this defines a sequence $\{u_n\}_n$ with $u_n \in \mathbb R^2$.
For $x=(x_1,x_2) \in \mathbb R^2$ we set $|x|^2= x_1^2+x_2^2$ and
the scalar product in $\mathbb R^2$ is $x\cdot y=x_1y_1+x_2y_2$.

Then, using the notations of Section 2.1, we define the basic space
$H$ as
\[
 H= \{u=(u_1, u_2, \ldots) \in (\mathbb R^2)^\infty:
 \sum_{n=1}^\infty |u_n|^2<\infty \}.
\]

The basis  in $H$ in given by the sequence
$\{e_1^{(1)}, e_1^{(2)},e_2^{(1)}, e_2^{(2)},e_3^{(1)}, e_3^{(2)},\ldots\}$
of elements of $(\mathbb R^2)^\infty$, where
\[
e_n^{(1)}=\left((0,0),\ldots,(0,0),(1,0),(0,0),\ldots \right)
\]
\[
e_n^{(2)}=\left((0,0),\ldots,(0,0),(0,1),(0,0),\ldots \right)
\]
with the nonvanishing vectors in place $n$.\\
The eigenvalues are
\[
 \lambda_n= k_0^2 \lambda^{2n}
\]

with $\lambda>1$. Hence we can take any $\alpha< \beta$ to fulfil
{\bf (C3)}. Inequality \eqref{semigr} holds with $c_{p,\nu}=(\frac p{e\nu})^p$.

We set $k_n=\sqrt \lambda_n$.
The bilinear term $B$ is defined by means of the components
$B_{n}=(B_{n,1}, B_{n,2})$ as follows (see, e.g., \cite{BF}):
\begin{align}
\begin{split}
&B_{1,1}(u,v)=ak_{2} [-u_{2,2}  v_{3,1} +u_{2,1}  v_{3,2}]
\\
&B_{1,2}(u,v)=-a k_{2} u_{2} \cdot v_{3}
\end{split}
&\\
\begin{split}
&B_{2,1}(u,v)=ak_{3} [-u_{3,2}  v_{4,1} +u_{3,1}  v_{4,2}]
           +b k_{2} [-u_{1,2} v_{3,1}+u_{1,1} v_{3,2}]
\\
&B_{2,2}(u,v)=-a k_{3} u_{3}\cdot  v_{4} -b k_{2}  u_{1} \cdot v_{3}
\end{split}
&\\
\intertext{and for $n >2$}
\begin{split}
&B_{n,1}(u,v)= \; ak_{n+1} [-u_{n+1,2}  v_{n+2,1} +u_{n+1,1}  v_{n+2,2}]
\\&\qquad\qquad\qquad   +b k_{n} [-u_{n-1,2} v_{n+1,1}+u_{n-1,1} v_{n+1,2}]
\\&\qquad\qquad\qquad           +a k_{n-1}[u_{n-1,2} v_{n-2,1}+u_{n-1,1} v_{n-2,2}]
\\&\qquad\qquad\qquad           +b k_{n-1}[u_{n-2,2} v_{n-1,1}+u_{n-2,1} v_{n-1,2}],
\end{split}
&\\
\begin{split}
&B_{n,2}(u,v)=-a k_{n+1} [u_{n+1,1}  v_{n+2,1} +u_{n+1,2}  v_{n+2,2}]
\\&\qquad\qquad\qquad    -b k_{n}  [u_{n-1,1} v_{n+1,1}+u_{n-1,2} v_{n+1,2}]
\\&\qquad\qquad\qquad    -a k_{n-1} [u_{n-1,1} v_{n-2,1}-u_{n-1,2} v_{n-2,2}]
\\&\qquad\qquad\qquad    -b k_{n-1} [u_{n-2,1} v_{n-1,1}-u_{n-2,2} v_{n-1,2}].
\end{split}&
\end{align}
where $a$ and $b$ are real numbers such that
\begin{equation}\label{a-b}
 a+b\lambda^{2\beta}=(a+b)\lambda^{4\beta}
\end{equation}
for some $\beta>0$, that is
\begin{equation}\label{condizbeta}
  \lambda^{2\beta}=-\frac a{a+b}
\end{equation}
(recall that $\lambda>1$). This condition implies
{\bf (C2 iii)}, whereas {\bf (C2 ii)}
holds for any real $a$ and $b$.
For instance, let us check that \eqref{a-b} implies {\bf (C2 iii)}.
We have
\[\begin{split}
  \sum_{n=1}^\infty k_n^{2\beta}& B_n(u,u) \cdot u_n
\\&
=\sum_{n=1}^\infty k_n^{2\beta} [B_{n,1}(u,u) u_{n,1}+B_{n,2}(u,u) u_{n,2}]
\\&
=\sum_{n=1}^\infty [a+b \lambda^{2\beta}-(a+b)\lambda^{4\beta}]
                \lambda k_n^{2\beta+1}
                (u_{n+2}\cdot u_n) (u_{n+1,2} + u_{n+1,1}) .
\end{split}\]

Moreover we have (see \cite{BF})
  \begin{lemma}
For any $\alpha_1, \alpha_2, \alpha_3 \in \mathbb R$
$$
 B:H^{\alpha_1}\times H^{\alpha_2} \to H^{-\alpha_3} \;
 \text{ with } \alpha_1+\alpha_2+\alpha_3\ge 1
$$
and there exists a constant $c$ (depending on $a,b,\lambda$ and  the
$\alpha_j$'s) such that
\[
 \|B(u,v)\|_{-\alpha_3}\le c \|u\|_{\alpha_1} \|v\|_{\alpha_2}
 \qquad  \forall u \in H^{\alpha_1}, v \in H^{\alpha_2}.
\]
\end{lemma}
This implies that conditions {\bf (C4)} are true:
\eqref{Bgenerale}  for any $\frac \beta 2
\le \alpha<\beta$ and \eqref{Ba} for any $\alpha$.

Condition \eqref{condE} holds for $\beta>\frac 12$; this includes the
interesting physical case of $\beta=1$ (see Section \ref{S-eqs}). Indeed,
for the SABRA shell model
\[
 B^{M}_{n,1}(x,x)-B_{n,1}(x,x)=
 \begin{cases}
 0 & \text{ for } n \le M-2\\
 -a k_M (x_{M,1}x_{M+1,2} - x_{M,2}x_{M+1,1}) & \text{ for } n=M-1\\
 -a k_{M+1}(x_{M+1,1}x_{M+2,2} - x_{M+1,2}x_{M+2,1}) &\\[-2mm]
      & \text{ for } n=M\\[-2mm]
     \;\qquad  -b k_M (x_{M-1,1}x_{M+1,2} - x_{M-1,2}x_{M+1,1}) &
 \end{cases}
\]
and
\[
 B^{M}_{n,2}(x,x)-B_{n,2}(x,x)=
 \begin{cases}
 0 & \text{ for } n \le M-2\\
 -a k_M (-x_{M,1}x_{M+1,1} - x_{M,2}x_{M+1,2})  & \text{ for } n=M-1\\
 -a k_{M+1}(-x_{M+1,1}x_{M+2,1} - x_{M+1,2}x_{M+2,2})&\\[-2mm]
      & \text{ for } n=M\\[-2mm]
  \;\qquad -b k_M (-x_{M-1,1}x_{M+1,1} - x_{M-1,2}x_{M+1,2}) &
 \end{cases}
\]
Therefore
\[
 \sum_{n=1}^M |B_n^{M}-B_n|^2=
  |B_{M-1}^{M}-B_{M-1}|^2+|B_{M}^{M}-B_{M}|^2
\]
so
\[
 \lim_{M\to \infty} \int \sum_{n=1}^M |B_n^{M}-B_n|^2 d\ms \le
  \lim_{M\to \infty}\frac{8}{\nu^2}
    [\frac{a^2}{\lambda^{2\beta}} k_M^{2-4\beta}
     + \frac{a^2}{\lambda^{2\beta}} k_{M+1}^{2-4\beta}
     + b^2 k_M^{2-4\beta}] =0.
\]
This holds for $\beta>\frac 12$.

We finally point out that our results of Section \ref{Sminv}
hold also in any space $\mathcal L^p(\ms)$
with $  p=1,2,\ldots$ (see Remark \ref{oss-Lp}).
Indeed, we have
\begin{equation}\label{sommabilita-p}
 \int |B_n(x,x)|^q \ \ms(dx)<\infty\qquad \forall n, q \in \mathbb N.
\end{equation}

\section{Inviscid models}

We are interested in the deterministic inviscid and  unforced dynamics
represented by  equation \eqref{eul}.
Here we present our results for the SABRA shell model with $\beta=1$
 (the physical relevant case)  only to make
simpler the exposition, but it can be generalized to the other fluid dynamic
models.

Equation \eqref{eul} is
formally obtained from equation \eqref{sequu} setting $\nu=0$ and
considering a vanishing right hand side.
More generally we can consider the nonlinear viscous equation
\begin{equation}\label{eq-ep}
 du^\ep(t)+[\nu \ep A u^\ep(t) + B(u^\ep(t),u^\ep(t))]dt
=  \sqrt {2\ep } \ dw(t), \qquad t>0.
\end{equation}
with $\ep>0$. When $\ep=0$ we get equation \eqref{eul} (with $\beta=1$).
Our results of the previuous sections hold true for
any $\ep>0$.

The fact that the measure $\msu$ is an invariant
measure for any $\ep>0$ can be easily checked.
We proceed as in the previous section, but now the
Kolmogorov operator associated to equation \eqref{eq-ep} is
$K^\ep=\ep Q+L$; bearing in mind \eqref{invQ} and \eqref{invL} we get that
$\msu$ is an infinitesimal invariant measure for the operator
$(K^\ep,FC^\infty_b)$. And for any $\ep>0$ the operator $(K^\ep,FC^\infty_b)$
is dissipative.

We are going to prove that  when the initial data is a random variable
 with law $\msu$, then equation \eqref{eul}
has a solution which is a stationary  random process, whose law at any
fixed time is $\msu$.

An important property is the integrability
of $B$ with respect to the measure $\msu$.
\begin{proposition}\label{Bstaz}
If $\nu>0$, then for any $\alpha<1$ we have
$$
 \int  \|B(x,x)\|^p_{\alpha} \ \msu(dx)<\infty
$$
for any $p\in \mathbb N$.
\end{proposition}

\proof We write the proof for $p=2$ but it is the same for the other
values of $p$, since $\msu$ is Gaussian and the $B_n$'s are
second order polynomial.
We have
\[\begin{split}
\int |B_{n,1}(x,x)|&^2 \msu(dx)
=
 \int |a k_{n+1} [-x_{n+1,2}  x_{n+2,1} +x_{n+1,1}  x_{n+2,2}]
\\&           +b k_{n} [-x_{n-1,2} x_{n+1,1}+x_{n-1,1} x_{n+1,2}]
\\&           +(a+b) k_{n-1}[x_{n-1,2} x_{n-2,1}+x_{n-1,1} x_{n-2,2}]|^2 \msu(dx)
\\&
\le 2 \int \{a^2 k^2_{n+1} [x^2_{n+1,2}  x^2_{n+2,1} +x^2_{n+1,1}  x^2_{n+2,2}]
\\&       +b^2 k^2_{n} [x^2_{n-1,2} x^2_{n+1,1}+x^2_{n-1,1} x^2_{n+1,2}]
\\&       +(a+b)^2 k^2_{n-1}[x^2_{n-1,2} x^2_{n-2,1}+x^2_{n-1,1} x^2_{n-2,2}]\}\msu(dx)
\\& =\frac{16}{\nu^2} \{a^2 k^2_{n+1}(\lambda_{n+1}\lambda_{n+2})^{-1}
                       +b^2 k^2_{n}(\lambda_{n-1}\lambda_{n+1})^{-1}
                       +(a+b)^2 k^2_{n-1}(\lambda_{n-1}\lambda_{n-2})^{-1}\}
\\&= \frac{4}{\nu^2k_0^2}
     \{a^2\lambda^{-4}+b^2+(a+b)^2\lambda^{4}\}\lambda^{-2n}.
\end{split}
\]

Similarly we estimate $\int |B_{n,2}(x,x)|^2 \mu^{1,\nu}(dx)$.
Therefore
\[\begin{split}
 \int \|B(x,x)\|_{\alpha}^2\msu(dx)
 &=
 \int \sum_{n=1}^\infty \lambda_n^{\alpha}|B_n(x,x)|^2 \msu(dx)\\
 &\le c_{\nu,k_0,\lambda} (|a|^2+|b|^2) \sum_{n=1}^\infty \lambda^{2n(\alpha-1)}
\end{split}
\]
which is finite if $\alpha<1$. \hfill $\Box$

Here is our main result.
\begin{theorem}
For any $\nu>0$,
there exists a $\msu$-stationary process, whose paths solve
equation \eqref{eul} $\po$-a.s.
In particular, the paths are in
$C^\delta(\mathbb R; H^{\alpha})$ (for any $0\le \delta<\frac 12$ and
$\alpha<1$).
\end{theorem}
\proof
We fix $\nu>0$ arbitrarily.
According to Corollary \ref{e-staz},
equation \eqref{eq-ep}
has a unique $\msu$-stationary solution $\overline v^{\ep}$;
this process is a strong solution
and has paths in $C([0,\infty);H^{\alpha})$  a.s..
(for $\alpha<1$,  but we always think of $\alpha$ as much close to 1
as possible).

First, we prove that the sequence $\{\overline v^{\ep}\}_{0<\ep\le 1}$
is tight in $C^{\tilde\delta}([0,T];H^{\tilde\alpha})$ for
any $\tilde\delta \in
(0,\frac 12)$ and $\tilde\alpha<\alpha$.

We write equation \eqref{eq-ep} in the mild form:
\begin{equation}\label{leqint}
 \ove(t)=\overline z^\ep(t)
         -\int_0^t e^{-\nu \ep A (t-s)} B(\ove(s), \ove(s))ds,
\end{equation}
where
\[
\overline z^\ep(t) =
e^{-\nu \ep A t} \ove(0)+\int_0^t e^{-\nu \ep A (t-s)} \sqrt {2\ep} dw(s)
\]
is the $\msu$-stationary solution of the linear equation
\[
  dz^\ep(t)+\nu \ep A z^\ep(t) dt
=  \sqrt {2\ep} \ dw(t)
\]
with the initial data of law $\msu$.

We consider the two terms in the right hand side of \eqref{leqint}.
Using the  $\msu$-stationarity  we have
that for any $0\le \delta < \frac 12$
there exists a constant $\overline C_\delta>0$ such that
\begin{equation}
 \sup_{0<\ep\le1}
  \mathbb E[\|\overline z^\ep\|_{C^\delta([0,T];H^{\alpha})}]
 \le  \overline C_\delta.
\end{equation}

We take $\eta\in (0,1)$ and set $\gamma=\alpha-2\eta$.
For the convolution integral in \eqref{leqint} we have

\begin{equation}
\begin{split}
&\Big\|\int_0^\cdot  e^{-\nu \ep A (\cdot-s)}
  B(\ove(s),\ove(s)) ds \Big\|^p_{W^{1,p}(0,T;H^\gamma)}\\
&=
\int_0^T \Big\|\int_0^t  e^{-\nu \ep A (t-s)}
  B(\ove(s),\ove(s)) ds \Big\|^p_{\gamma} dt
 +\int_0^T \Big\|B(\ove(t),\ove(t))\Big\|^{p}_{\gamma}dt\\
&+ \int_0^T \Big\| \int_0^t \nu \ep A e^{-\nu \ep A (t-s)}
  B(\ove(s),\ove(s)) ds \Big\|^p_{\gamma} dt
\\
& \le
\int_0^T t^{p-1} \left(\int_0^t
  \| e^{-\nu \ep A (t-s)}B(\ove(s),\ove(s))\|^p_{\gamma} ds\right)  dt
 +\int_0^T \Big\|B(\ove(t),\ove(t))\Big\|^{p}_{\gamma}dt\\
&+ \nu \ep
 \int_0^T \left(\int_0^t
 \|A e^{-\nu \ep A(t-s)}B(\ove(s),\ove(s))\|_\gamma ds\right)^p dt
\\
&\le
 \int_0^T t^{p-1} \left(\int_0^t  \| B(\ove(s),\ove(s))\|^p_{\gamma} ds\right)  dt
 +\int_0^T \Big\|B(\ove(t),\ove(t))\Big\|^{p}_{\gamma}dt\\
&+\nu \ep \int_0^T
 \left(\int_0^t \|A^{1-\eta} e^{-\nu \ep A(t-s)}A^\eta
 B(\ove(s),\ove(s))\|_\gamma ds\right)^p dt
\\
& \le
 (\frac 1p T^p +1) \int_0^T \Big\|B(\ove(t),\ove(t))\Big\|^{p}_{\gamma}dt
 + \nu \ep \int_0^T
\left(\int_0^t c_{p,\nu}
  \frac{\|B(\ove(s),\ove(s))\|_\alpha}{(t-s)^{1-\eta}} ds\right)^p dt
  \text{ by } \eqref{semigr}.
\end{split}\end{equation}
For the latter integral we use H\"older inequality and get that
\[
\left(\int_0^t
  \frac{\|B(\ove(s),\ove(s))\|_\alpha}{(t-s)^{1-\eta}} ds\right)^p
\le \left(\int_0^t \frac{ds}{(t-s)^{1-\frac \eta 2}}\right)^{2p\frac {1-\eta}{2-\eta}}
\left(\int_0^t\|B(\ove(s),\ove(s))\|_\alpha
   ^{\frac 2\eta -1}ds\right)^{p\frac \eta{2-\eta}}.
\]
Hence, for $p>\frac 2 \eta-1$ we have
\begin{multline}\label{nu-Balpha}
\Big\|\int_0^\cdot  e^{-\nu \ep A (\cdot-s)}
  B(\ove(s),\ove(s)) ds \Big\|^p_{W^{1,p}(0,T;H^\gamma)}
\\\le
 (\frac 1p T^p +1) \int_0^T
 \Big\|B(\ove(t),\ove(t))\Big\|^{p}_{\gamma}dt
 +\nu \ep T^m \int_0^T \|B(\ove(t),\ove(t))\|^p_\alpha dt
\end{multline}
for some positive constant $m=m_{\eta,\nu,p}$.

Integrating with respect to the measure $\ms$ and using the invariance
we get
\begin{equation}\label{nu-Balpha2}
\begin{split}
&\mathbb E \Big\|\int_0^\cdot  e^{-\nu \ep A (\cdot-s)}
  B(\ove(s),\ove(s)) ds \Big\|^p_{W^{1,p}(0,T;H^{\gamma})}\\
& \qquad\le T(1+\frac 1p T^p+\nu \ep T^m) \int \|B(x,x)\|_{\alpha}^p \msu(dx)
\end{split}
\end{equation}

Now, we use that $W^{1,p}(0,T)\subset C^\delta([0,T])$ if  $1-\frac 1p>\delta$.
Then, using the previous estimates in  \eqref{leqint},
given any $0\le\delta<\frac 12$,
$p>\frac 1{1-\delta}$ and $p>\frac 2 \eta-1$ we have
\begin{equation}\label{tigh}
  \sup_{0 < \ep \le 1}
  \mathbb E[\|\ove\|^p_{C^\delta([0,T];H^\gamma)}]<\infty .
\end{equation}
On the other hand, the space $C^\delta([0,T];H^\gamma)$ is
compactly embedded in $C^{\tilde \delta}([0,T];H^{\tilde \gamma})$
if $\tilde \delta <\delta$ and $\tilde \gamma<\gamma$; this follows from the
compact embedding $H^\gamma \Subset H^{\tilde \gamma}$
and from the Ascoli-Arzel\`a theorem.
Because these results hold for any $\delta \in
[0,\frac 12)$ and $\tilde \gamma<\gamma<\alpha<1$
(with $p$ big enough, but we use \eqref{sommabilita-p}),
we can consider any $\tilde \delta
  <\frac 12$ and
 any $\tilde \gamma<1$.
The tightness  follows from \eqref{tigh} as usual  by means of
Chebyshev inequality. And to simplify notation henceforth we consider the
tightness
in the space $C^\delta([0,T];H^\alpha)$
($\delta <\frac 12$ and $\alpha<1$).

By the tightness result and Prohorov theorem,
the sequence of the laws of $\ove$  has a subsequence
$\{\overline v^{\ep_n}\}_{n=1}^\infty$ weakly
convergent as $n \to \infty$ (with $\ep_n \to 0$)
in $C^\delta([0,T];H^{\alpha})$ to some limit measure.
By a diagonal argument, this holds for any $T$ and therefore the limit
measure  leaves in  $C^\delta([0,\infty); H^{\alpha})$.
By Skorohod theorem, there
exist a probability space
$(\tilde \Omega, \tilde {\mathbb F}, \tilde \po)$,
a random variable $\tilde v$ and a sequence $\{\tilde v^{\ep}\}$
such that law($\tilde v^{\ep} $)=law($\overline v^{\ep}$),
law($\tilde v $)=$\msu$ and $\tilde v^{\ep}$ converges
to $\tilde v$ a.s. in
$C^\delta([0,\infty); H^{\alpha})$.

We now identify the
equation satisfied by $\tilde v$. We are going to prove that
  $\tilde \po$-almost
each path solves \eqref{eul}.

It is enough to control the behavior of the terms with $B$. First
\[\begin{split}
e^{-\nu \ep A (t-s)}B(\tilde v^{\nu,\ep}(s),&\tilde v^{\nu,\ep}(s))
   -B(\tilde v^{\nu}(s),\tilde v^{\nu}(s))
\\&=
e^{-\nu \ep A (t-s)}
\big[B(\tilde v^{\nu,\ep}(s),\tilde v^{\nu,\ep}(s))
     -B(\tilde v^{\nu}(s),\tilde v^{\nu}(s))\big]
\\&+
 \big[ e^{-\nu \ep A (t-s)}-I\big]
 B(\tilde v^{\nu}(s),\tilde v^{\nu}(s)).
\end{split}\]
When we consider the second addend in the mild form expression,
it trivially converges to zero; but for the convergence
of the first one
it is enough to verify that
\[
\int_0^t \|B(\tilde v^{\nu,\ep}(s),\tilde v^{\nu,\ep}(s))
     -B(\tilde v^{\nu}(s),\tilde v^{\nu}(s))\|_{\alpha-1} ds \to 0
\]
as $\ep \to 0$; for this we use
the bilinearity and the estimate \eqref{alpha-alpha}.

Similarly we work on the time interval $[-T,0]$ by considering
the reversed-time parabolic nonlinear equation
\begin{equation}\label{eq--ep}
 du^\ep(t)+[-\nu \ep A u^\ep(t) + B(u^\ep(t),u^\ep(t))]dt
=  \sqrt {2\ep} \ dw(t),\qquad t<0
\end{equation}
It has a unique $\msu$-stationary solution $\underline v^{\ep}$;
this process is a strong solution,
has paths in $C^\delta((-\infty,0];H^{\alpha})$.
The tigthness and the convergence are obtained in the same way as above.
\hfill $\Box$\\
\\
\textbf{Acknowledgments:} The work of H. Bessaih was partially
supported by the NSF grant No. DMS 0608494.

\end{document}